\documentclass[letterpaper, 10 pt, conference]{ieeeconf}  % Comment this line out
\IEEEoverridecommandlockouts                              % This command is only
\usepackage{graphics} % for pdf, bitmapped graphics files
\usepackage{epsfig} % for postscript graphics files
\usepackage{mathptmx} % assumes a new font selection scheme installed
\usepackage{times} % assumes a new font selection scheme installed
\usepackage{amsmath} % assumes amsmath package installed
\usepackage{amssymb}  % assumes amsmath package installed
\usepackage{todonotes}
\usepackage{booktabs}
\usepackage[linesnumbered,ruled,vlined]{algorithm2e}
\usepackage{multirow}
\usepackage{hyperref}

\SetKwInput{KwInput}{Input}                % Set the Input
\SetKwInput{KwOutput}{Output}              % set the Output

\title{\LARGE \bf
Addressing Discrete Dynamic Optimization via a Logic-Based Discrete-Steepest Descent Algorithm}

\author{Zedong Peng, Albert Lee, David E. {Bernal Neira}% <-this % stops a space
\thanks{Zedong Peng, Albert Lee, and David E. {Bernal Neira} are with the Davidson School of Chemical Engineering, Purdue University, West Lafayette, IN, 47907 USA {\tt\small \{peng372,lee4382,dbernaln\}@purdue.edu}}%
}

\begin{document}

\maketitle
\thispagestyle{empty}
\pagestyle{empty}

%%%%%%%%%%%%%%%%%%%%%%%%%%%%%%%%%%%%%%%%%%%%%%%%%%%%%%%%%%%%%%%%%%%%%%%%%%%%%%%%
\begin{abstract}

Dynamic optimization problems involving discrete decisions have several applications, yet lead to challenging optimization problems that must be addressed efficiently.
Combining discrete variables with potentially nonlinear constraints stemming from dynamics within an optimization model results in mathematical programs for which off-the-shelf techniques might be insufficient.
This work uses a novel approach, the Logic-based Discrete-Steepest Descent Algorithm (LD-SDA), to solve Discrete Dynamic Optimization problems. The problems are formulated using Boolean variables that enforce differential systems of constraints and encode logic constraints that the optimization problem needs to satisfy.
By posing the problem as a generalized disjunctive program with dynamic equations within the disjunctions, the LD-SDA takes advantage of the problem's inherent structure to efficiently explore the combinatorial space of the Boolean variables and selectively include relevant differential equations to mitigate the computational complexity inherent in dynamic optimization scenarios.
We rigorously evaluate the LD-SDA with benchmark problems from the literature that include dynamic transitioning modes and find it to outperform traditional methods, i.e., mixed-integer nonlinear and generalized disjunctive programming solvers, in terms of efficiency and capability to handle dynamic scenarios. 
This work presents a systematic method and provides an open-source software implementation to address these discrete dynamic optimization problems by harnessing the information within its logical-differential structure.
% Implementing large-scale dynamic simulation is becoming more common in process systems, showing that more advanced optimization methods are being adopted and used effectively.
%For example, the optimization problem is formulated into mixed-integer dynamic optimization (MIDO), which can tackle various dynamic problems with diverse numerical solutions.
% Due to the intricate nature of MIDO problems, Logic MIDO (LMIDO) formulations have emerged, utilizing logic variables to streamline and diminish their computational complexity.
%There are many ways to solve the LMIDO problems.
%However, solving the LMIDO problem is an open challenge.
%Here, we show the Logic-based Discrete-Steepest Descent Algorithm (LD-SDA) to tackle the LMIDO problem.
%We benchmarked the existing problem from the previous literature and implemented the problem with other solution methods, including LD-SDA.
%Compared to other existing methods, the LD-SDA performed more efficiently than the different approaches.
%Moreover, we formulated more complex problems than the previous ones and solved them with various methods.
%LD-SDA solved the Logic-based Dynamic Optimization problem more productively than other methods.
% Two or three more sentences on the broader perspectives

\end{abstract}

%%%%%%%%%%%%%%%%%%%%%%%%%%%%%%%%%%%%%%%%%%%%%%%%%%%%%%%%%%%%%%%%%%%%%%%%%%%%%%%%

\section{Introduction}
Control tasks with high-level discrete or logical decisions, such as integrated process design and control~\cite{flores2007simultaneous}, trajectory optimization~\cite{landry2016aggressive}, and energy management in hybrid electric vehicles~\cite{robuschi2021multiphase}, can be modeled using Mixed-Integer Dynamic Optimization (MIDO).
% Large-scale dynamic simulation has emerged as an essential tool for optimizing chemical process systems, supporting the growing acceptance and application of advanced optimization techniques in this domain.
% Central to these advancements is Mixed-Integer Dynamic Optimization (MIDO), a sophisticated framework designed to tackle the complexities inherent in systems that encompass discrete decisions and continuous variables. 
% MIDO provides a robust mechanism for enhancing the static design of process systems and their dynamic operational control~\cite{richards2005mixed}.
% This ability to handle planning and real-time control makes MIDO an essential tool for optimizing process controls.
% The scope of MIDO extends across a diverse array of process tasks, effectively addressing challenges in start-up and shut-down procedures~\cite{raghunathan2004mpec}, grade transitions~\cite{prata2008integrated}, batch operations, and more. 
% MIDO enables seamless coordination between batch scheduling and control systems. 
% It also allows for the simultaneous planning and management of processes~\cite{flores2006simultaneous}.
% These capabilities highlight MIDO's wide-ranging role in enhancing chemical processing efficiency, reliability, and effectiveness.
However, MIDO problems are computationally challenging to solve, and standard optimization solvers cannot simultaneously handle discrete variables and differential algebraic equations (DAE). 
% Even though MIDO is a powerful tool for handling both discrete decisions and differential algebraic equations (DAE), solving MIDO problems is computationally challenging, especially when the ordinary differential equations (ODE) are nonlinear.
% The ability to handle both discrete decisions and differential algebraic equations (DAE) makes MIDO an essential tool for mixed-integer control problems.
% Since the discrete behaviors in mixed-integer control problems can sometimes be modeled as logic variables and constraints, an alternative to MIDO is Differential-Algebraic Generalized Disjunctive Programming (DAGDP), also known as logic Mixed-Integer Dynamic Optimization (LMIDO).
% Several numerical methods have been devised to practically address it~\cite{flores2007simultaneous,bansal2002case,biegler2010nonlinear}.
The standard routine for solving MIDO problems is to discretize differential equations, resulting in a mixed-integer nonlinear programming (MINLP) problem~\cite{flores2007simultaneous,bansal2002case,biegler2010nonlinear}.
However, the resulting MINLP problem is usually nonconvex, as it involves nonlinear equality constraints, and existing MINLP solvers are insufficient to solve these problems~\cite{kronqvist2019review}, particularly in online settings relevant to MIDO.
% Traditionally, the discrete in mixed-integer control problems often involved models that used logic variables as decision variables.
% This approach has given rise to logic MIDO (LMIDO) problems.
% LMIDO problems emerge when the optimization model incorporates logical decisions typically represented through disjunctions.

Since the discrete behaviors in mixed-integer control problems can sometimes be modeled as logic variables and constraints, an alternative to MIDO is Differential-Algebraic Generalized Disjunctive Programming (DAGDP).
In DAGDP problems, logical decisions, represented as Boolean variables indicating the choices via True or False values, interact with the dynamic systems depicted by differential equations.
Compared to the MIDO formulation, the disjunctions in DAGDP problems allow for a more compact representation of the logic-induced control problem, where certain operations or configurations are only feasible or relevant under specific discrete conditions.
Moreover, another benefit of the DAGDP formulation is that it facilitates the use and development of alternative solution strategies in addition to the MINLP reformulation, which more effectively leverages the structure of the optimization model.
% However, solving the LMIDO problem efficiently remains an open challenge.
The use of dynamic models within the framework of DAGDP offers a promising avenue to model control problems associated with processes, and this paper introduces an innovative approach to effectively tackling their complexities.

% LMIDO can be solved by converting into MIDO and transcribing the problem into MINLP problems.
% The MINLP can be solved by applying outer approximation, Benders decomposition, or branch and bound.
% Another approach to solving the LMIDO problem is to discretize it. 
% LMIDO problems are converted into nonlinear programming problems that are less complex and do not use discrete decision variables.
% Though using complementarity constraints introduces non-convexities that increase the complexity of the NLP problem, this approach remains promising and warrants further research efforts.

% In dynamic optimization, MINLP is employed to handle problems involving switching variables and capturing mode transitions of the system.
% On the other hand, GDP focuses only on relevant constraints and reduces the problem's complexity compared to MINLP, avoiding redundancy.
% This can make GDP more advantageous than MINLP for complex dynamic systems.

\subsection{Related Work}
\label{sec:related work}
The approaches to solving DAGDP problems primarily derive from and combine established ideas from both generalized disjunctive programming (GDP) and dynamic optimization.
Typically, solution strategies for dynamic optimization problems include single shooting, multiple shooting, direct transcription, and numerical integration~\cite{biegler2010nonlinear}.
In the literature, direct transcription, e.g., finite difference or orthogonal collocation, is first applied to discretize the DAE system.
Through this, the DAGDP problem will be reformulated into a GDP model, which will be solved by reformulation into MINLP (BigM or Hull) or logic-based methods~\cite{flores2007simultaneous,ruiz2014logic}. 
% For example, Ruiz-Femenia et al. first apply orthogonal collocation to discretize the DAEs and then use a logic-based outer approximation method to solve LMIDO problems~\cite{ruiz2014logic}. 
% Flores-Tlacuahuac and Biegler develop three MINLP reformulations of the LMIDO problem, i.e., BigM, Hull, and binary multiplication reformulations~\cite{flores2007simultaneous}. The effectiveness of these reformulations is further evaluated through the outer approximation and NLP branch-and-bound algorithms.
Although this approach is simple and intuitive, it has a significant drawback, that the reformulated MINLP problem is usually nonconvex.

In addition to MINLP reformulation, an alternative method to solve MIDO problems is to apply complete discretization techniques based on complementarity constraints, also known as the Mathematical Program with Equilibrium Constraints (MPEC)~\cite{baumrucker2008mpec}.
The complementarity constraint avoids the use of discrete variables.
However, MPEC reformulations yield a highly nonconvex nonlinear programming (NLP) problem, usually challenging to solve.

Other DAGDP solution approaches in the literature handle the discrete variables first and then solve a series of continuous dynamic optimization subproblems.
This strategy does not necessarily depend on a particular discretization technique and provides more freedom to tackle the DAE system.
For example, the MIDO problems can be decomposed into a primal dynamic optimization and mixed-integer linear programming (MILP) sub-problems.
These can be solved later using a MINLP decomposition algorithm, such as the outer-approximation and generalized bender decomposition methods~\cite{burnak2020mixed}.
Chachuat et al.~\cite{chachuat2006global} address MIDO problems by adapting the outer-approximation algorithm and applying the branch-and-bound method to solve dynamic optimization subproblems to global optimality.
Deriving valid linear inequalities to construct the MILP subproblems is a challenging task for MIDO problems, which usually have nonconvex nonlinear constraints.
Inspired by this approach, we propose decomposition methods that decouple the dynamic optimization from the discrete choices. We aim to tackle the discrete variable optimization without defining a MILP.

Another technique used to optimize discrete problems is through a discrete steepest descent algorithm (D-SDA), where the search space is mapped to a lattice representing each realization of the discrete variables and its corresponding subproblem over the remaining continuous decision, followed by a series of steepest descent steps to find improving solutions~\cite{linan2020optimalI}.
This lattice can be smaller than the original combinatorially large one defined by all combinations. 
The D-SDA is designed to address MINLP problems and has been used to efficiently solve dynamic optimization problems that integrate design and nonlinear model predictive control (NMPC)~\cite{linan2022discrete} and consider cases even under uncertainty~\cite{palma2023simultaneous}. 
This approach effectively avoids the pitfalls of nonconvexities and suboptimal solutions inherent in other optimization methods.
In this particular case, the dynamic model is not directly enforced by the discrete choices in the problem.
While the D-SDA demonstrates robustness in managing complex decision variables and nonconvex problems, as illustrated in a distillation column case study, the authors acknowledge the computational intensity of the subproblems as a challenge and point towards future enhancements to improve efficiency and applicability to broader systems.

Another recent extension of the D-SDA is to consider the case where discrete choices imply sets of (potentially nonlinear) constraints, usually expressed as GDP problems~\cite{bernal2022process}.
This approach was named logic-based D-SDA (LD-SDA).
Using the LD-SDA shows that leveraging the logical structure in these disjunctive programs leads to an improved solution of highly nonlinear problems, such as a catalytic distillation column, compared to using the D-SDA over a reformulated MINLP.
The constraints implied by the Boolean variables addressed by the LD-SDA were algebraic, and no previous work has addressed the case where a system of differential equations implies these constraints.

\subsection{Contributions}
This work's contributions can be summarized as follows:
\begin{itemize}
    \item We extend the LD-SDA to solve DAGDP problems. The LD-SDA utilizes a logic-based search strategy to explore the search space given by the discrete choices, here Boolean variables, involved in disjunctions, and it is independent of the type of method used to solve the dynamic optimization subproblems. This flexibility allows a variety of dynamic optimization approaches to be applied, such as finite difference methods, orthogonal collocation, and even numerical integrators.
    % (mention explicitly that we do not need to do orthogonal collocation and we could have had a numerical integrator solving the subproblems)
    \item This work integrated the model transformations in the open-source algebraic modeling language Pyomo~\cite{bynum2021pyomo}, Pyomo.DAE and Pyomo.GDP, for DAGDP modeling. Several numerical instances of mode transitionded.
    The implementation is available at \href{https://github.com/SECQUOIA/LD-SDA-Dynamic}{github.com/SECQUOIA/LD-SDA-Dynamic}.
    \item An open-source and general implementation of LD-SDA is provided as an option for the GDPOpt solver~\cite{chen2022pyomo} for users to solve GDP problems, potentially with DAE systems in the disjunctions.
\end{itemize}
% The novelty of this project is to show that the LD-SDA can solve the dynamic optimization problem.
% LD-SDA solved the GDP problems that involved stationary models and proved that these problems can be solved more efficiently than the D-SDA when solving complex problems.
% We performed a case study on three-stage dynamic model switching adapted from the problem formulation presented in~\cite{ruiz2014logic}.
% Moreover, to find out the performance of the LD-SDA on more complicated dynamic optimization problems, we formulated a multiple-stage dynamic optimization problem with three selection modes and solved the problem with different solvers.

\section{Background}

\subsection{Discrete-Continuous Dynamic Optimization Problems}

In this work, we focus on DAGDP problems with Boolean variables, i.e., \(Y \in \{False, True \} = \{\bot, \top\}\), and differential-algebraic equations (DAE) in the disjunctions.
The problem formulation can be written as follows.
\begin{subequations}
\label{DAGDP}
\begin{align}
\min_{\substack{x(t), y(t), \\ u(t), p , \\ Y(t)\in \{\bot, \top\}}} \ & \psi(x(t), y(t), u(t), p, t) \label{DAGDP:obj} \\
\text{s.t.} \ & \xi(x(t),\dot{x}(t), y(t), u(t), p, t) = 0\label{DAGDP:constr1}\\
&\bigvee_{r \in D_k}
\begin{bmatrix}
    Y_{kr} \\
    \varphi_{kr} (x(t),\dot{x}(t), y(t), u(t), p, t) = 0
\end{bmatrix}\label{DAGDP:constr2}\\
&\veebar_{r \in D_k} (Y_{kr}) = True = \top \quad \forall k \in K \label{DAGDP:constr3}\\
&\Omega(Y_{kr}) = True = \top \label{DAGDP:constr4}\\
& x(0) = x_0 \label{DAGDP:constr5}\\
&x^L \leq x(t) \leq x^U , y^L \leq y(t) \leq y^U \label{DAGDP:constr6}\\
&u^L \leq u(t) \leq u^U, p^L \leq p \leq p^U. \label{DAGDP:constr7}
\end{align}
\end{subequations}
where \(x(t)\) represents the time-dependent state vectors, with time derivatives indicated by \(\dot{x}(t)\). 
\(y(t)\) denotes the algebraic variable states, and \(u(t)\) the control actions.
The vector \(p\) encompasses the decision variables of the system not dependent on time \(t\). 
Eq.~\eqref{DAGDP:constr1} represents the DAE system or the purely algebraic equations that are consistently enforced.
% The case where DAE systems are outside the disjunctions is not explicitly considered here. However, one can obtain the following formulation by discretizing those equations first as described in \S~\ref{sec:discretization}, hence obtaining Eq.~\eqref{DAGDP:constr1}.
For each disjunction \(k \in K\), a selection is made from a set of options defined by \(D_k\).
For every disjunct \(r \in D_k\), a Boolean variable \(Y_{kr}\), differential algebraic equations, representing system dynamics, are specified as Eq.~\eqref{DAGDP:constr2}.
% The case where DAE systems only exist outside the disjunctions is not explicitly considered here.
The objective function \(\psi(\cdot)\) and functions \(\xi (\cdot)\) and \(\varphi_{kr} (\cdot)\) that define constraints are potentially nonlinear.
The Boolean variables follow an exclusive OR \((\veebar)\) constraint as Eq.~\eqref{DAGDP:constr3}, ensuring that within each disjunction \(k\), exactly one option \(r\) is selected.
If \(Y_{k,r}(t)\) is set to \(True = \top\), the corresponding disjunct and the constraints inside will be active.
Otherwise, the disjunct can be discarded.
Eq.~\eqref{DAGDP:constr4} is the logical propositions that represent the relationships of the logical variables. The initial state and the bounds are provided in Eqs.~\eqref{DAGDP:constr5} - \eqref{DAGDP:constr7}.
% In the disjunction of Boolean variables, the GDP arranges these variables to reflect sequential choices, such as identifying distinct locations or selecting specific moments.
% This ordering also facilitates the quantification of task repetitions within the decision-making process.

% Orthogonal collocation in~\cite{biegler2010nonlinear}.
\subsection{Reformulations and logic-based methods for GDP}
\label{sec:GDPsol}

% Grossmann GDP review (include it briefly; see PSE paper~\cite{bernal2022process})
There are different strategies to solve GDP problems.
One approach is to reformulate the GDP problem into MINLP problems via BigM and Hull reformulations. 
The BigM reformulation introduces a sufficiently large scalar that makes the particular constraint redundant when its indicator variable is not selected. 
The Hull reformulation lifts the model to a higher-dimensional space by introducing copies of the continuous variables and constraints inside disjunctions. 
The Hull reformulation always yields tighter relaxations than the BigM reformulation at the expense of larger model sizes.
For more details, see~\cite{grossmann2021advanced}.

% M is a sufficiently large scalar that makes the particular constraint redundant when its indicator variable is not selected (i.e., yi = 0). 

% The two classic reformulations of the GDP are the BigM and the Hull reformulation. The BigM reformulation 
% In the BigM approach inequalities are relaxed by setting \(1-y\) times a large constant \(M\) that effectively renders them non-binding when \(y=0\).
% The Hull reformulation introduces copies of the continuous variables and constraints, which are activated through a perspective reformulation of the problem.

% The BigM reformulation requires fewer variables than the Hull reformulation, while the Hull reformulation has a tighter continuous relaxation than the BigM reformulation.
% Even though MINLP reformulations can be the default solution for GDP problems, these methods introduce numerous algebraic constraints, some of which might not be relevant to a particular solution and might even lead to numerical instabilities.

In addition to MINLP reformulation, one can exploit explicit logical propositions in GDP problems via logic-based methods, such as the logic-based outer approximation (LOA) and logic-based branch and bound (LBB)~\cite{chen2022pyomo}.
Contrary to the MINLP reformulation, logic-based methods formulate specific subproblems corresponding to the values of the logical variables while solving the problems.
These subproblems only include constraints activated by the logical variables within each evaluated assignment of $True/False$ to each Boolean variable, or configuration.
For instance, if the specific logical configuration \(\hat{Y}\) is given, the disjunctions can be fixed, and the subproblem becomes
\begin{subequations}
\label{GDP subproblem}
\begin{align}
&\min_{\substack{x(t), y(t), u(t), p}} \quad \psi(x(t), y(t), u(t), p, t) \label{GDP subproblem:obj} \\
&\text{s.t.} \quad \xi(x(t),\dot{x}(t), y(t), u(t), p, t) = 0 \label{GDP subproblem:constr1}\\
&\varphi_{kr} (x(t),\dot{x}(t), y(t), u(t), p, t) = 0 \ \forall Y_{kr} = \top \label{GDP subproblem:constr2}\\
& x(0) = x_0 \label{GDP subproblem:constr3}\\
&x^L \leq x(t) \leq x^U , y^L \leq y(t) \leq y^U \label{GDP subproblem:constr4}\\
&u^L \leq u(t) \leq u^U, p^L \leq p \leq p^U. \label{GDP subproblem:constr5}
\end{align}
\end{subequations}
% \begin{equation}
% \begin{align}
% \min_{\substack{x_{ij}, y_{ij}, u_{ij}}} \quad& Z = \psi (x_{ij}, y_{ij}, u_{ij}, p) \\
% \text{s.t.} \quad &F(x_{ij}, y_{ij}, u_{ij}, p) \leq 0 \\
% &H(x_{ij}, y_{ij}, u_{ij}, p) \leq 0 \\
% &G(x_{ij}, y_{ij}, u_{ij}, p) \leq 0 \\
% &f_{kr} (x_{ij}, y_{ij}, u_{ij}, p) = 0 \ \text{if}~ \hat{\mathbf{Y}}_{ij} = True \quad r \in D_k, k \in K \\
% &\varphi_{kr} (x_{ij}, y_{ij}, u_{ij}, p) = 0  \ \text{if}~ \hat{\mathbf{Y}}_{ij} = True \quad r \in D_k, k \in K \\
% &x_{ij} \in \mathbb{R}^n; \ Y_{kr} \in \{True, False\} \\ 
% &\ i \in N_e, \ j \in N_c, \ k \in K, \ r \in D_k. \\
% \end{align}
% \label{eq:sub}
% \end{equation}

The subproblem represents the optimization problem under constraints with a fixed logical configuration.
This problem might avoid evaluating numerically challenging nonlinear equations whenever their corresponding logical variables are irrelevant.
Since the subproblem satisfies the logical proposition \eqref{DAGDP:constr4}ithm avoids solving subproblems from infeasible logical configurations.

Solving the subproblems~\eqref{GDP subproblem}, which result from exploration of the space defined by discrete variables, can result in convergence to the optimal solution of \eqref{DAGDP}.
The methods for choosing the series of subproblems lead to different logic-based methods, among them LOA and LBB.
LOA uses gradient-based linearization of the nonlinear constraints at the optimal solution of Eq.~\eqref{GDP subproblem} to approximate the feasible region of the original problem.
The additional constraint would be added to a mixed-integer programming problem, denoted as the main problem.
The optimal solution to the main problem returns a configuration of Boolean variables.
On the other hand, LBB systematically solves GDP by exploring the values of Boolean variables in the search tree, where each node represents the partial fixation of the Boolean variables.
The solutions in the node provide bounds to the optimal solution.
Both methods are designed to efficiently find the optimal configuration of Boolean variables~\cite{chen2022pyomo}.

\subsection{Discretization for Dynamic Optimization}
\label{sec:discretization}

DAGDP problems contain differential and algebraic constraints in their disjunctions.
An approach to obtain a problem with only algebraic constraints that can be solved as a GDP, as described in \S~\ref{sec:GDPsol}, is to use the transcription approach~\cite{biegler2010nonlinear}, which transforms sets of differential equations into algebraic equations through an orthogonal collocation within finite elements. The transformed GDP becomes
\begin{subequations}
\label{GDP}
\begin{align}
&\min_{\substack{x_{ij}, y_{ij}, u_{ij}, p , \\ Y_{kr}\in \{\bot, \top\}}} \quad \psi(x_{ij}, y_{ij}, u_{ij}, p) \label{GDP:obj} \\
&\text{s.t.} \qquad F(x_{ij}, y_{ij}, u_{ij}, p) \le 0 \label{GDP:constr1}\\
&\bigvee_{r \in D_k}
\begin{bmatrix}
    Y_{kr}\\
    f_{kr} (x_{ij}, y_{ij}, u_{ij}, p) = 0
\end{bmatrix} \label{GDP:constr2}\\
&\veebar_{r \in D_k} (Y_{kr}) = True = \top \quad \forall k \in K\\
&\Omega(Y_{kr}) = True = \top \label{GDP:constr3}\\
& x(0) = x_0 \\
&x^L \leq x_{ij} \leq x^U , y^L \leq y_{ij} \leq y^U \\
&u^L \leq u_{ij} \leq u^U, p^L \leq p \leq p^U \\
& i \in N_e, \ j \in N_c, \ k \in K, \ r \in D_k,
\end{align}
\end{subequations}
where \(N_e\) is the number of finite elements and \(N_c\) is the number of internal collocation points used for properly discretizing the DAE.
When comparing the DAGDP formulation with GDP, the differential equations in Eqs.~\eqref{DAGDP:constr1} and \eqref{DAGDP:constr2} are mapped into algebraic equations in Eqs.~\eqref{GDP:constr1} and \eqref{GDP:constr2}.

When discretizing the DAGDP into GDP, the entire time horizon is separated into several finite elements.
The dynamic behavior of the process is captured in each stage using a series of points given by an orthogonal collocation.
The smoothness of the dynamic response determines the right number of these elements.
The size of a given finite element \(i\) represents the particular length of the independent variable.
The orthogonal collocation within each finite element facilitates the precise determination of the internal location points, ensuring accurate modeling of the process dynamics.

When applying both transformations, the subproblems become the following NLP problems
\begin{subequations}
\label{NLP}
\begin{align}
&\min_{\substack{x_{ij}, y_{ij}, u_{ij}, p}} \quad \psi(x_{ij}, y_{ij}, u_{ij}, p) \label{NLP:obj} \\
&\text{s.t.} \quad F(x_{ij}, y_{ij}, u_{ij}, p) \leq 0 \label{NLP:constr1}\\
&f_{kr} (x_{ij}, y_{ij}, u_{ij}, p) = 0 \quad \forall Y_{kr} = \top  \label{NLP:constr2}\\
% &\varphi_{ij,kr} (x_{ij}, y_{ij}, u_{ij}, p) = 0 \quad \forall Y_{ij,kr} = \top \label{NLP:constr3}\\
& x(0) = x_0 \\
& x^L \leq x_{ij} \leq x^U , y^L \leq y_{ij} \leq y^U \\
& u^L \leq u_{ij} \leq u^U, p^L \leq p \leq p^U \\
& i \in N_e, \ j \in N_c, \ k \in K, \ r \in D_k.
\end{align}
\end{subequations}

\subsection{Convergence Analysis}
The motivation for the LD-SDA is based on discrete convex analysis~\cite{murota1998discrete}, where the convergence to global optimal points of various convex functions (M-convex and L-convex) over lattices is characterized by local optimality over integral neighborhoods.
The value of the function corresponds to the optimal solution of a dynamic optimization problem given a fixed set of discrete choices.
However, we cannot guarantee that the solutions to the DAGDP models over the discrete choices define a function that is either M- or L-convex.
Therefore, the LD-SDA serves as a heuristic method, and the convergence of global optimality cannot be guaranteed.

% The direct approach, first discretize then optimize, does not explicitly guarantee the uniqueness of global optima. 
% However, this method provides a framework for achieving globally optimal solutions within a small margin of error. 
% Thus, discretizing and optimizing effectively address global optimality concerns within the practical limitations of the problem and computational methods employed
%~\cite{sager2015efficient}.

\section{Logic-Based Discrete-Steepest Descent Algorithm}

In this section, we explain how to apply the Logic-based Discrete-Steepest Descent Algorithm (LD-SDA) to tackle DAGDP problems with ordered Boolean variables.

As mentioned in \S\ref{sec:related work}, the key of LD-SDA is to first map the ordered Boolean variables into a lattice with each point representing a particular realization of the disjunctions. 
Then, LD-SDA performs the neighbor search and the line search to find improving solutions by solving the corresponding continuous dynamic optimization subproblems.

\subsection{Reformulation}
In the DAGDP model, ordered Boolean variables \( \{{Y}_{kr} \mid r \in D_k\} \) can be reformulated into a set of integer variables referred to as \textit{external variables}.
% The reformulation is also possible across other ordered sets of Boolean variables.
For example, in each exclusive OR constraint \eqref{DAGDP:constr3}, the Boolean variables \(Y_{1k}, Y_{2k}, \dots, Y_{\lvert D_k \rvert k}\)  can be represented by discrete values \(z_k \in \{1,2,\dots,\lvert D_k \rvert \}\) according to an ordered sequence they follow. 
After the reformulation, the feasible region of the boolean variables is mapped into a \(\lvert D_k \rvert\)-dimensional integer lattice.
% Suppose there is a Boolean variable \(Y_{ak}\) where \(a\) is the element of set \(S = \{1,2,\dots,\lvert D_k \rvert \}\) in disjunction \(k\).
% The Boolean variable \(Y_{ak}\) is reformulated into an external variable \(\mathbf{z}_E = a\).
% The Boolean variables reformulation allows us to map the Boolean variables into a lattice of external variables.

\subsection{Algorithm Description}
In the lattice, each point corresponds to a continuous dynamic optimization subproblem by fixing the disjunctions, which can be solved by arbitrary dynamic optimization methods.
The LD-SDA starts from the given initial point and obtains the initial primal bound (PB) by solving the subproblem.
% At the initial point, the disjunctions in the DAGDP model can be fixed and the resulting dynamic optimization subproblem can be solved via, say, discretization methods. 
Then, a neighbor search over the external variable lattice is performed to find the steepest descent direction.
Two types of integral neighborhoods, i.e., search directions, are supported, defined by the $L_2$ and $L_\infty$ norms.
If no better solution is found in the neighbor search, a locally optimal solution is reached, and the algorithm terminates.
Otherwise, there exists at least one improving direction. In this case, the algorithm will move to the best neighbor, and a line search will be performed in the improving direction.
If a worse solution is detected during the line search, a new neighbor search and line search will be repeated at the incumbent best-found point until a local optimum is found.
After each neighbor and line search, the explored points will be added to the set $G$ of explored points, used to avoid exploring the same point twice.
The detailed steps are described in Algorithm \ref{alg:D-SDA}. For simplicity, we use $\mathbf{z}$ to denote the external variables and $\mathbf{x}$ for all remaining variables in the DAGDP problem.

\begin{algorithm}[htb]
\caption{Logic-based Discrete-Steepest Descent Algorithm (LD-SDA) for DAGDP problems}
\label{alg:D-SDA}
\DontPrintSemicolon
  
  \KwInput{ An external variable feasible solution $\mathbf{{z}_{0}}$; Integral neighborhood $\in \{L_2, L_\infty \}$.}
  Initialize: $k \leftarrow 0$, $G \leftarrow \{\mathbf{{z}_0} \}$ \\
  Generate the search directions $d \in D$.\\
  Solve the initial DO subproblem with the given $\mathbf{z_0}$ 
  $(\mathbf{x_{0}} , \text{PB}_0) \leftarrow \textbf{SolveDO}(\mathbf{z_0})$ \\
  \While {$True$}
  {
  Generate neighbors $N_k = \{n: n=\mathbf{{z}_k} + d \ \forall d \in D\} \setminus G$ \\
  % Find $D_k(\mathbf{{z}_E}) = \{\mathbf{d} : \; \mathbf{\alpha} - \mathbf{{z}_E} = \mathbf{d}, \forall \ \mathbf{\alpha} \in N_k(\mathbf{{z}_E})\}$ \\
  Perform \textbf{Neighbor Search} 
  $(\mathbf{z_{k+1}}, \mathbf{x_{k+1}}, \mathbf{d^*}, G, \text{PB}_{k+1}) \leftarrow \textbf{NS}(\mathbf{{z}_k}, N_k, \text{PB}_k)$\\
  $k \leftarrow k+1$\\
  % \If{\normalfont{\texttt{PB improved during NS}}}
  \If{$\text{PB}_{k} > \text{PB}_{k-1}$}
    {
         % Set \texttt{lineSearch} $\leftarrow True$ \\
         \While {$\text{PB}_{k} > \text{PB}_{k-1}$}
            {
                Perform \textbf{Line Search}
                $(\mathbf{z_{k+1}},\mathbf{x_{k+1}}, G,\text{PB}_{k+1}) \leftarrow \textbf{LS}(\mathbf{{z}_{k}}, \mathbf{d^*}, \text{PB}_k)$\\
                $k \leftarrow k+1$\\
            }
    }
    \Else
    {
    $\mathbf{z^*} \leftarrow \mathbf{z_{k}}$, 
      $\mathbf{x^*} \leftarrow \mathbf{x_{k}}$\\
    \Return Best found feasible solution $\mathbf{z^*}$ and $\mathbf{x^*}$
      
    }
  }
   % \KwOutput{Best found feasible solution $\mathbf{z^*}$ and $\mathbf{x^*}$}
\end{algorithm}

\section{Computational experiments}

\subsection{Implementation Details}
This section shows the effectiveness of the LD-SDA through several computational experiments. The DAGDP models are written using the DAE and GDP modules in Pyomo. 
We provide an implementation of LD-SDA in GDPOpt in Pyomo. 
We further benchmark LD-SDA against solving the problem as a reformulated MINLP problem, and other logic-based methods, such as LOA, GLOA, and logic-based enumeration. 
KNITRO and BARON are used as (MI)NLP solvers. 
For KNITRO, \texttt{mip\_multistart} is set to 1 to enable a mixed-integer multi-start heuristic and improve the chances of finding the global solution. 
All tests ran on a Linux cluster with 48 AMD EPYC 7643 2.3GHz CPUs and 1 TB RAM, restricted to using only a single thread. 
% The code for the algorithm and models is available at \href{https://github.com/SECQUOIA/LD-SDA-Dynamic}{https://github.com/SECQUOIA/LD-SDA-Dynamic}.

\subsection{Three-stage Dynamic Model Switching}
Consider the optimization of a system with two dynamic modes and three stages~\cite{ruiz2014logic}. 
At each stage \(s\), only one of the two modes can be enforced. 
This problem aims to compute the dynamic model and the optimal control actions that apply in each stage by maximizing the square of the state variable over the time horizon, and is formulated as

\begin{subequations}
\label{ex1}
\begin{align}
    &\min_{\substack{x(t), u(t), \\ Y \in \{\bot, \top\}}} V(x) = - \int_{t_0}^{t_s} x^2(t) \label{ex1:obj}\\
    &\text{s.t.} \begin{bmatrix}
        Y_{s,1} \\
        \frac{dx}{dt} = -x e^{x-1} + u
    \end{bmatrix} \veebar
    \begin{bmatrix}
        Y_{s,2} \\
        \frac{dx}{dt} = \frac{0.5 x^3 + u}{20}
    \end{bmatrix}, t \in  [t_{s-1},t_s] \\
        & Y_{s,2} \Rightarrow \lor_{s'<s} Y_{s',1} \label{eq:sequence1-1}\\
    & Y_{s,2} \Rightarrow \lnot \lor_{s'>s} Y_{s',1} \label{eq:sequence1-2}\\
    & t_0 = 0, \ t_s = s \qquad \qquad \qquad \forall s = \{1, 2, 3\} \\
    & x(0) = 1, \ u(0) = 4,
\end{align}
\end{subequations}
where $x$ is the state variable and $u$ is the control variable. 
To solve this problem, we first use orthogonal collocation to discretize the differential equations in each disjunction with 30 finite elements and three collocation points at each stage. 
The discretized model has 553 variables, 819 constraints, and three disjunctions.
To apply the LD-SDA to this problem, we first reformulate the disjunctions using external variable $z_{s} \in \{1,2\}$. 
$z_{s} = 1$ represents when mode 1 is active and $z_{s} = 2$ when mode 2 is active at stage $s$. 

The computational results for \eqref{ex1} are presented in Table \ref{tab:three stage results ordering}, and the time limit is set at 900 seconds. 
All the methods using BARON reach the maximum time limit and cannot prove global optimality within the time limit. KNITRO fails to solve the BigM and Hull reformulation and returns the infeasible termination condition, while all logic-based methods using KNITRO find the optimal solution within 5 seconds. 
In this problem, both LD-SDA $L_2$ and $L_\infty$ explore all feasible disjunctions similarly to an enumeration over the logic space. 
LOA and GLOA terminate after around 5 seconds, slower than logic-based enumeration and D-SDA.
These computational results show the advantage of the logic-based method over MINLP reformulations for this problem.

% \begin{table}[]
%     \caption{Computational results of problem \eqref{ex1:obj} (without ordering constraint)}
%     \label{tab:three stage results}
%     \centering
%     \begin{tabular}{llrrl}
%     \toprule
%     Strategy & Solver & Obj & Time & Status \\
%     \midrule
%     MINLP BigM & KNITRO & - & 0.04 & Infeasible \\
%     MINLP BigM & BARON & -8.08 & 900.13 & maxTimeLimit \\
%     MINLP Hull & KNITRO & - & 2.93 & Infeasible \\
%     MINLP Hull & BARON & -9.91 & 900.37 & maxTimeLimit \\
%     L-Enumerate & KNITRO & -12.74 & 3.09 & Optimal \\
%     L-Enumerate & BARON & -9.18 & 900.22 & maxTimeLimit \\
%     LOA & KNITRO & -12.74 & 11.57 & Optimal \\
%     LOA & BARON & -9.18 & 900.75 & maxTimeLimit \\
%     GLOA & KNITRO & -12.74 & 12.45 & Optimal \\
%     GLOA & BARON & -9.18 & 900.74 & maxTimeLimit \\
%     LD-SDA $L_2$ & KNITRO & -12.74 & 2.78 & Optimal \\
%     LD-SDA $L_2$ & BARON & -9.18 & 900.22 & maxTimeLimit \\
%     LD-SDA $L_\infty$ & KNITRO & -12.74 & 3.21 & Optimal \\
%     LD-SDA $L_\infty$ & BARON & -9.18 & 900.33 & maxTimeLimit \\
%     \bottomrule
%     \end{tabular}
% \end{table}

\begin{table}[]
    \centering
    \caption{Computational results of problem \eqref{ex1}}
    \label{tab:three stage results ordering}
    \begin{tabular}{llrrl}
    \toprule
    Strategy & Solver & Obj & Time [s] & Status \\
    \midrule
    \multirow{2}{*}{MINLP BigM} & KNITRO & - & 0.05 & Infeasible \\
     & BARON & -9.18 & 900+ & maxTimeLimit \\
    % \midrule
    \hline
    \multirow{2}{*}{MINLP Hull} & KNITRO & - & 9.75 & Infeasible \\
     & BARON & -12.74 & 900+ & maxTimeLimit \\
    % \midrule
    \hline
    \multirow{2}{*}{L-Enumerate} & KNITRO & -12.74 & 1.56 & Optimal \\
     & BARON & -9.18 & 900+ & maxTimeLimit \\
    % \midrule
    \hline
    \multirow{2}{*}{LOA} & KNITRO & -12.74 & 4.63 & Optimal \\
     & BARON & -12.74 & 900+ & maxTimeLimit \\
    % \midrule
    \hline
    \multirow{2}{*}{GLOA} & KNITRO & -12.74 & 5.03 & Optimal \\
     & BARON & -12.74 & 900+ & maxTimeLimit \\
    % \midrule
    \hline
    \multirow{2}{*}{LD-SDA $L_2$} & KNITRO & -12.74 & 1.50 & Optimal \\
     & BARON & -9.18 & 900+ & maxTimeLimit \\
    % \midrule
    \hline
    \multirow{2}{*}{LD-SDA $L_\infty$} & KNITRO & -12.74 & 1.47 & Optimal \\
     & BARON & -9.18 & 900+ & maxTimeLimit \\
    \bottomrule
    \end{tabular}
\end{table}

\subsection{Multi-stage Dynamic Model Switching}
Consider the following DAGDP with three dynamic modes and $S$ stages and sequencing constraints:
\begin{subequations}
\label{ex2}
\begin{align}
    &\min_{\substack{x(t), u(t), \\ Y \in \{\bot, \top\}}} V(x) = - \int_{t_0}^{t_S} x^2(t) \tag{6a} \label{ex2:obj} \,\\
    &\text{s.t.} 
    \begin{bmatrix}
        Y_{s,1} \\
        \frac{dx}{dt} = \frac{-x}{e^{1-x}} + u
    \end{bmatrix} \veebar 
    \begin{bmatrix}
        Y_{s,2} \\
        \frac{dx}{dt} = \frac{0.5 x^3 + u}{20}
    \end{bmatrix} \veebar
    \begin{bmatrix}
        Y_{s,3} \\
        \frac{dx}{dt} = \frac{x^2 + u}{t + 20}
    \end{bmatrix}, \notag\\
    & t \in[t_{s-1}, t_s]\tag{6b} \\
    & Y_{s,2} \Rightarrow \lor_{s'<s} Y_{s',1} \label{eq:sequence1} \tag{6c}\\
    & Y_{s,2} \Rightarrow \lnot \lor_{s'>s} Y_{s',1} \label{eq:sequence2} \tag{6d} \\
    & Y_{s,3} \Rightarrow \lor_{s'<s} Y_{s',2} \label{eq:sequence3} \tag{6e} \\
    & Y_{s,3} \Rightarrow \lnot \lor_{s'>s} Y_{s',2} \label{eq:sequence4} \tag{6f} \\
    & t_0 = 0, t_s = s \qquad \qquad \qquad \forall s = \{1, ..., S\} \tag{6g} \\
    & x(t_0) = 1, u_{t} \in [-4, 4], x \in [0, 10]. \tag{6h}
\end{align}
\end{subequations}
Eqs.~\eqref{eq:sequence1}-\eqref{eq:sequence4} are mode sequence constraints which denote that mode 1 should be performed before mode 2 and mode 2 before mode 3.
We discretize the differential equations with 30 finite elements and three collocation points at each stage using orthogonal collocation. 
The statistics of the discretized models are presented in Table \ref{tab:statistic of the multi-stage dynamic problem}.

\begin{table}[hb]
    \centering
    \caption{Statistic of the multi-stage dynamic switching problem \eqref{ex2}}
    \label{tab:statistic of the multi-stage dynamic problem}
    \begin{tabular}{cccc}
        \toprule
         \# of stages & \# of variables & \# of constraints & \# of disjunctions \\
         \midrule
         4 & 741 & 1183 & 4 \\
         5 & 927 & 1547 & 5 \\
         6 & 1113 & 1911 & 6 \\
         7 & 1299 & 2275 & 7 \\
         8 & 1485 & 2639 & 8 \\
         9 & 1671 & 3003 & 9 \\
         \toprule
    \end{tabular}
\end{table}

\begin{figure}
    \centering
    \includegraphics[width=0.46\textwidth]{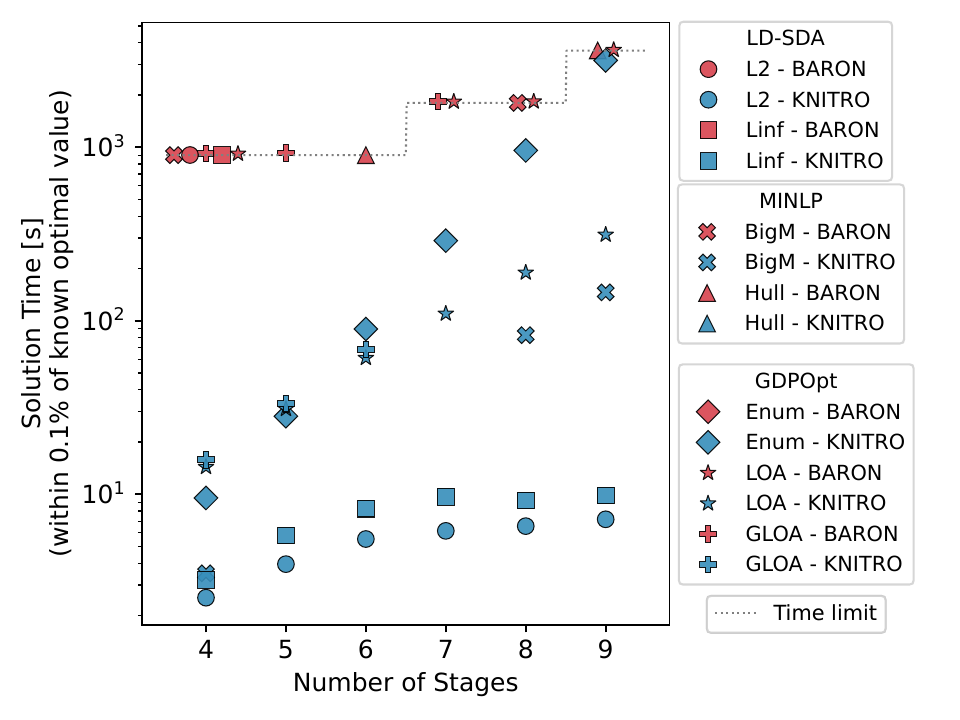}
    \caption{Computational results of MINLP reformulation, LD-SDA, and other logic-based methods for varying \(S\) in problem \eqref{ex2}}
    \label{fig:computational results}
\end{figure}

Using LD-SDA to solve this problem requires a reformulation based on external integer variables. 
A straightforward reformulation is to define a variable $z_{s} \in \{1,2,3\}$ to represent the mode choice at each stage. 
However, this would result in a search space as high-dimensional as the original Boolean variables space. 
This space is limited to three or fewer points in each dimension, limiting the possibility of line search and resulting in slow convergence depending on the initial point. 
Therefore, we propose another reformulation based on the mode transition, which can also apply to other problems with sequence constraints, such as in scheduling and planning. 
Instead of using the integer variable to denote the mode choice at each stage, we use it to indicate whether the mode transition happens at each stage. 
Due to the sequence constraint, only two transitions are allowed. i.e., $\{A: 1 \to B: 2 \to 3\}$, defining external variables $z_{s,A}, z_{s,B}$. 
Compared to the first reformulation, the high-dimensional space is reduced into a two-dimensional space. 
The relationship between mode transition and selection is expressed in Eqs.~\eqref{eq:mode transition1}-\eqref{eq:mode transition2}. Moreover, since mode transition might not occur within the stage horizon, we define extra variables in Eqs.~\eqref{eq:mode transition3}-\eqref{eq:mode transition4} to capture this case and force one of the $z_{s}$ variables of each transition to be $True$ in Eqs.~\eqref{eq:mode transition5}-\eqref{eq:mode transition6}.
\begin{align}
    z_{s, A} \Leftrightarrow Y_{s-1,1} \wedge Y_{s,2} & \quad \forall s=\{2, ..., S\} \tag{6i} \label{eq:mode transition1}\\
    z_{s, B} \Leftrightarrow Y_{s-1,2} \wedge Y_{s,3} & \quad \forall s=\{3, ..., S\}\tag{6j}  \label{eq:mode transition2}\\
    z_{S+1, A} \Leftrightarrow \lnot \vee_{s \in \{2,\dots,S\}}  {Y_{s,2}} & \tag{6k} \label{eq:mode transition3}\\
    z_{S+1, B} \Leftrightarrow \lnot \vee_{s \in \{3,\dots,S\}}  {Y_{s,3}} & \tag{6l} \label{eq:mode transition4}\\
    \veebar_{s \in \{2,\dots,S+1\}} z_{s, A} = True = \top & \tag{6m} \label{eq:mode transition5}\\
    \veebar_{s \in \{3,\dots,S+1\}} z_{s, B} = True = \top & \tag{6n} \label{eq:mode transition6}
\end{align}
The computational results of solving problem \eqref{ex2} with $S=\{4,\dots,9\}$ using LD-SDA, MINLP reformulation, and other logic-based methods are presented in Fig.~\ref{fig:computational results}. 
For the 4-stage problems, all the methods using BARON, except logic-based enumeration, can find the optimal solution. 
However, with increasing number of stages, none of these methods could guarantee to find the global optimal solution by closing the optimality gap within the time limit of up to 1 hour. 
Both LD-SDA $L_2$ and $L_{\infty}$ using KNITRO as NLP solver converge to the optimal solution for all problems and outperform the other logic-based and MINLP reformulation methods. 
% With the increase of stages, the solution time of the logic-based enumeration method and the LOA method increases exponentially.

\begin{figure}
    \centering
    \includegraphics[width=0.48\textwidth]{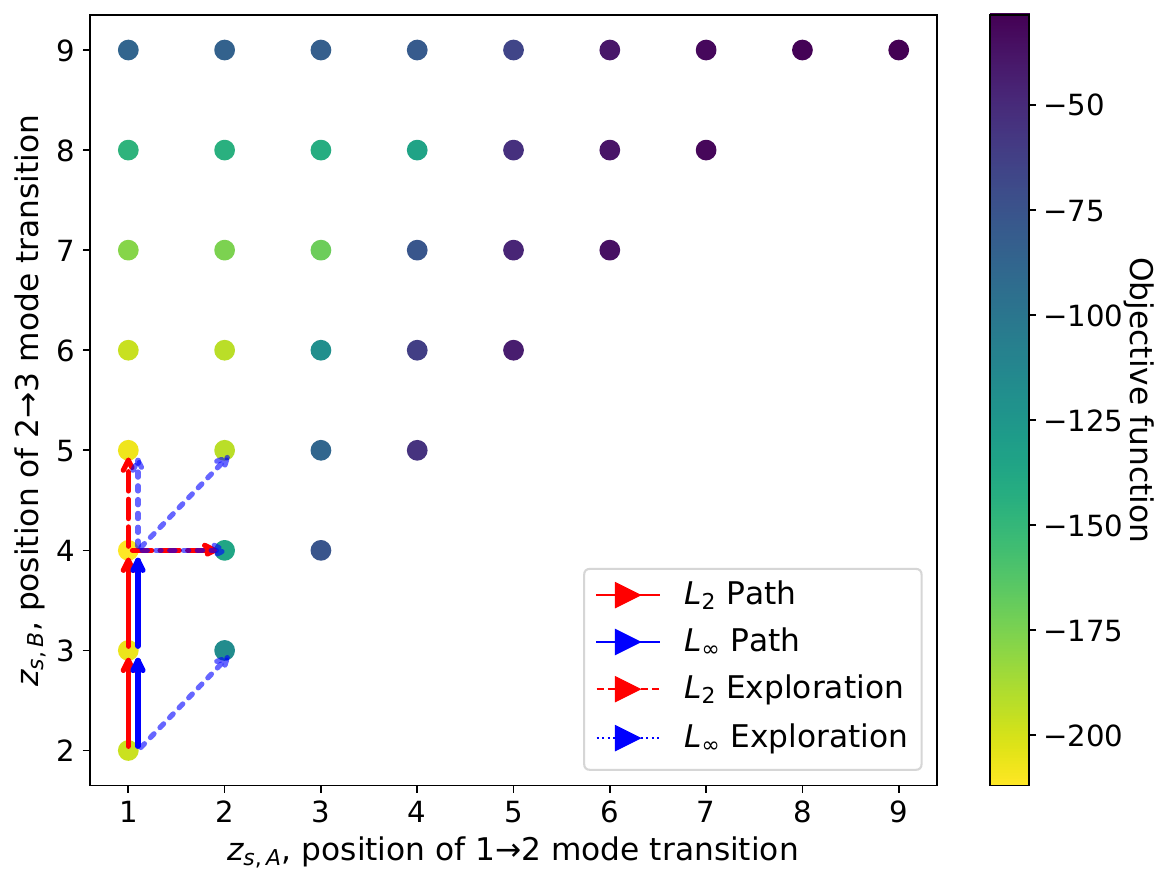}
    \caption{Search path of the LD-SDA for problem \eqref{ex2} with \(S=9\)}
    \label{fig:LDSDA L2 and Linf search path}
\end{figure}

The search path of LD-SDA $L_2$ and $L_{\infty}$ for the nine-stage problem~\eqref{ex2} is presented in Fig.~\ref{fig:LDSDA L2 and Linf search path}. 
Both algorithms start from point \((1,2)\), which means that the $A:1 \to 2$ mode transition occurs at stage 2, and the $B: 2 \to 3$ mode transition occurs at stage 3. 
After one round of neighbor search and a following line search, both LD-SDA $L_2$ and $L_{\infty}$ converge to the optimal point \((1, 4)\). 
Since LD-SDA $L_{\infty}$ allows more search directions, it explores two more points than LD-SDA $L_2$ before termination. 
However, neither method converging to the global optimum is guaranteed.
If the optimal point is not in the search directions of the starting or intermediate points, the LD-SDA could be trapped in a local optimum.

\section{Conclusions and Future Work}

% \subsection{Conclusions}

This work presents the LD-SDA as an optimization tool for discrete dynamic optimization modeled via Differential Algebraic Generalized Disjunctive Programming (DAGDP). 
The LD-SDA exploits the logical structure of the problem and allows more flexibility in handling the DAE systems. 
Two examples based on dynamic mode selection are provided to evaluate the performance of LD-SDA. 
The computational results demonstrate that LD-SDA outperforms the equivalent MINLP formulation and other logic-based methods when DAEs are discretized via direct transcription.
Moreover, an open-source implementation of the LD-SDA is provided, and this work integrates the model transformations in Pyomo.DAE and Pyomo.GDP for DAGDP modeling.

Future research directions include the integration of the LD-SDA with decomposition-based methods and exploring the effects of cutting planes on the LD-SDA.

%%%%%%%%%%%%%%%%%%%%%%%%%%%%%%%%%%%%%%%%%%%%%%%%%%%%%%%%%%%%%%%%%%%%%%%%%%%%%%%%
% \section{ACKNOWLEDGMENTS}

% The authors gratefully acknowledge the contribution of the National Research Organization and the reviewers' comments.

%%%%%%%%%%%%%%%%%%%%%%%%%%%%%%%%%%%%%%%%%%%%%%%%%%%%%%%%%%%%%%%%%%%%%%%%%%%%%%%%
\bibliographystyle{IEEEtran}
\bibliography{references}
%%%%%%%%%%%%%%%%%%%%%%%%%%%%%%%%%%%%%%%%%%%%%%%%%%%%%%%%%%%%%%%%%%%%%%%%%%%%%%%%

\end{document}